\documentclass[reqno, 11pt]{amsart}
\usepackage{graphicx, enumerate, enumitem, amsmath, amssymb,amsthm,mathabx,manfnt,hyperref,braket}

\usepackage[OT2,T1]{fontenc}
\DeclareSymbolFont{cyrletters}{OT2}{wncyr}{m}{n}
\DeclareMathSymbol{\Sha}{\mathalpha}{cyrletters}{"58}

\newtheorem{thm}{Theorem}[section]

\newtheorem{prop}[thm]{Proposition}

\theoremstyle{definition}

\theoremstyle{remark}

\numberwithin{equation}{section}

\newcommand{\abs}[1]{\left\vert#1\right\vert}
\newcommand{\pair}[2]{\left\langle#1,#2\right\rangle}
\newcommand{\Dd}{\mathcal{D}}

\newcommand{\cx}{{\mathbb{C}}}

\newcommand{\Mm}{\mathfrak{M}}

\newcommand{\e}{\epsilon}

\newcommand{\dbar}{\overline{\partial}}

\newcommand{\bc}{{\mathsf{bc}}}

\newcommand{\dist}{{\rm dist}}

\title[Distributional Boundary Values]{Distributional boundary values : some new perspectives}\thanks{This work was partially supported by a grant from the Simons Foundation (\#316632 to Debraj Chakrabarti). Debraj Chakrabarti was also partially supported by an Early Career internal grant from Central Michigan University. Rasul Shafikov was partially supported by an NSERC grant.}
\author{Debraj Chakrabarti}
\address{Department of Mathematics, Central Michigan University, Mt. Pleasant,  MI 48859,  USA}
\email{chakr2d@cmich.edu}
\author{Rasul Shafikov}
\address{Department of Mathematics, University of Western Ontario,  London, Ontario, Canada, N6A 5B7 }
\email{shafikov@uwo.ca}
\begin{document}
\maketitle

\section{Boundary values of holomorphic functions as currents}
Given a domain in a complex space, it is a fundamental problem to identify the class of boundary values of holomorphic functions on the domain. This notion is widely used in complex analysis, from the Cauchy integral formula to characterization of boundaries of complex subvarieties (Harvey-Lawson \cite{hala1}). For smoothly bounded domains in $\mathbb C^n$, boundary values are usually understood as a subclass of the so-called CR functions on the 
boundary, i.e., those satisfying the tangential Cauchy-Riemann equations. If the boundary is of class 
$\mathcal C^\infty$, then one may consider the Cauchy-Riemann equations in the weak sense, which gives rise to CR distributions. It is known (cf. \cite{straube84}) that  for a bounded domain with $\mathcal{C}^\infty$  boundary in $\mathbb C^n$, $n>1$, 
every holomorphic function of polynomial growth (i.e., the growth of the function near the boundary is bounded by some power of the distance to the boundary) admits a  boundary value which is a CR distribution. Distributional boundary values on  generic  CR submanifolds of higher codimension exist also for holomorphic functions of polynomial growth defined on a wedge attached to the submanifold,  (see \cite{baouendibook}).
There is also a parallel theory of generalized functions, the Sato hyperfunctions, which allows one to consider boundary values of arbitrary holomorphic functions on domains with real-analytic boundaries (cf. \cite{polkingwells}).

 It is natural to ask whether a notion of generalized boundary values of holomorphic functions exists for domains with 
 nonsmooth boundary. At the outset it is clear that as we reduce the regularity of the boundary, the class of holomorphic functions which admit boundary values would also become smaller.  In \cite{chak-sh} we define boundary values as $(0,1)$-currents in the ambient manifold satisfying certain conditions. This approach allows us to define boundary values on domains not necessarily with smooth boundary, in particular prove the existence of boundary values on domains with generic corners.
To formulate this result, assume that $\Omega$ is a relatively compact  domain in a complex  manifold  
$\Mm$ given in the form $\Omega = \bigcap_{j=1}^N \Omega_j$,
where each $\Omega_j\subset\Mm$ is a  smoothly bounded domain.  If  for each subset $S\subset\{1,\dots,N\}$ the intersection $B_S=\bigcap_{j\in S} b\Omega_j$, if non-empty,
is a CR manifold of CR-dimension $n-\abs{S}$, we say that  $\Omega$ is  a {\em domain with generic corners.}
The primary example of  domains with generic corners  are product domains. We denote by ${\rm dist}(z,X)$ the distance from a point $z\in\Mm$ to a set $X$ induced by some metric on $\Mm$ compatible with its topology.
We say that  a holomorphic $f\in \mathcal{O}(\Omega)$ is of {\em polynomial growth} if there is a $C>0$ and $k\geq 0$ such that we have for each $z\in \Omega$ that
\[ \abs{f(z)}\leq \frac{C}{\dist(z,\partial\Omega)^k}.\]
We denote the space of holomorphic functions of polynomial growth on $\Omega$ by $\mathcal{A}^{-\infty}(\Omega)$. 

\begin{thm}\label{thm-bcexistence} Let $\Omega$ be a domain with generic corners in a complex manifold $\Mm$ of complex dimension $n$, and let $f\in \mathcal{A}^{-\infty}(\Omega)$. There is a $(0,1)$-current $\bc f \in \Dd'_{0,1}(\Mm)$ such that the following 
holds. If $U$ is a coordinate neighbourhood of $\Mm$, and $\psi\in\Dd^{n,n-1}(\Mm)$  is a smooth $(n,n-1)$ form which has support
in $U$, and there is a vector $v\in \cx^n$ such that in the coordinates on $U$, the vector $v$ points outward from $\Omega$ along each $\partial\Omega_j$  inside $U$,
 then we have
 \begin{equation}\label{eq-bc}
 \pair{\bc f}{\psi}= \lim_{\epsilon\downarrow 0}\int_{\partial\Omega} f_\epsilon \psi ,
\end{equation}
where   $f_\epsilon(z)= f(z-\epsilon v)$.
\end{thm}

For proof and further discussion, see \cite{chak-sh}. It is also shown there that provided $\bc f$ exists, it is unique.  We refer to $\bc f$ as the {\em boundary current} induced by the holomorphic function $f$ of polynomial growth.
It is immediate from the formula \eqref{eq-bc} that for holomorphic functions that extend continuously to $\partial \Omega$ we
simply have $\bc f = f[\partial\Omega]^{0,1}$, where $[\partial\Omega]$ is the 1-current of integration on $\partial\Omega$, i.e., $\pair{[\partial\Omega]}{\phi}= \int_{\partial\Omega}\phi$ for a smooth $(2n-1)$-form $\phi$ of compact support, and for a 1-current $\gamma$, we denote by $\gamma^{0,1}$ the 
$(0,1)$-part of this current. One can also see by a use of Stokes' formula that for a holomorphic function $f$ on $\Omega$, which belongs to $L^1(\Omega)$
(with respect to any Riemannian measure on $\Mm$), we have $\bc f = -\dbar (f[\Omega])$, where $[\Omega]$ is the 0-current of integration on $\Omega$.
This even makes sense when $\Omega$ is an arbitrary open relatively compact subset of $\Mm$.

  It natural therefore to ask whether in Theorem~\ref{thm-bcexistence} the condition of generic corners on the domain $\Omega$ is necessary or not.
  It turns out that  if we want all holomorphic functions 
of polynomial growth on $\Omega$ to have boundary currents, then the condition that the boundary of $\Omega$ has generic corners is necessary, and 
the proof of this fact is the main result of this note:
\begin{prop} There is a complex manifold $\Mm$, a piecewise smooth domain (with non-generic corners) $\Omega\Subset\Mm$ and a holomorphic function
$f$ of polynomial growth on $\Omega$, such that $\bc f$ does not exist.
\end{prop}
In a later note, we will show that much more is true: on each piecewise smooth domain with non-generic corners, there is a holomorphic function of polynomial growth which does not admit a boundary current.

\begin{proof} Let $\Mm=\cx$ and 
\[
\Omega = \{ x+iy \in \mathbb C : |x-1|<1, |y-1|<1\}.\]
We will show that the function $f(z)=z^{-2}$, which is holomorphic in $\Omega$ and is of polynomial
growth there, does not admit the boundary value current as defined in Theorem~\ref{thm-bcexistence}.  Suppose to 
the contrary that $\bc f$ exists. Let $U= \{\abs{z}< \frac{3}{2}\}$. The vector $v=-(1+i)$ points outward from $\Omega$ along $\partial\Omega\cap U$, and therefore, for each $\psi\in \mathcal{D}^{1,0}(U)$, we have
\[ 
\pair{\bc f}{\psi} = \lim_{\epsilon\downarrow 0}\int_{\partial\Omega} f_\epsilon \psi,
\]
where $f_\epsilon(z)= f(z-\epsilon v)$. We choose $\psi$ to be $x\,dz$ in a neighbourhood of  the closed unit disc $\{\abs{z}\leq 1\}$ and vanishing outside $U$. We will show that 
\begin{equation}\label{nolimit}
\lim_{\e\downarrow 0}\int_{\partial \Omega} \frac{\psi}{(z-\e v)^2}
\end{equation}
does not exist, this will disprove the existence of $\bc f$.  Writing
\[ \int_{\partial \Omega} \frac{\psi}{(z-\e v)^2}  = \int_{\partial \Omega\cap\{\abs{z}\leq 1\}} \frac{\psi}{(z-\e v)^2} + 
\int_{\partial \Omega\cap\{\abs{z}> 1\}} \frac{\psi}{(z-\e v)^2},\]
we note that the second integral remains bounded as $\epsilon\to 0$, so it suffices to show that the first integral goes to infinity as $\epsilon\to 0$.
We have,
\begin{align*}
\int_{\partial \Omega\cap\{\abs{z}\leq 1\}} \frac{\psi}{(z-\e v)^2}&  = \int_{\partial \Omega\cap\{\abs{z}\leq 1\}} 
\frac{x(dx+idy)}{(z+\e+i\e)^2}\\ &= \int_0^1 \frac{xdx}{(x+\e+i\e)^2}\\
& =
\int_0^1 \frac{x^2(x+2\e)dx}{((x+\e)^2+\e^2)^2}
-2i\e \int_0^1 \frac{(x^2+\e x)dx}{((x+\e)^2+\e^2)^2} .
\end{align*}
Consider the real part of the last line, which we write as
$$
I+II = \int_0^1 \frac{x^3dx}{((x+\e)^2+\e^2)^2} +
2\e \int_0^1 \frac{x^2 dx}{((x+\e)^2+\e^2)^2}.
$$
Direct computation shows that 
$$
\int \frac{x^3dx}{((x+\e)^2+\e^2)^2}  = \frac{1}{2}\ln (x^2+2x\e+2\e^2)-2\tan^{-1} 
\left(\frac{x+\e}{\e}\right)  +\frac{\e x}{x^2+2x\e+2\e^2}+C.
$$
Therefore,
\[ I= \frac{1}{2} \ln \left( \frac{1}{2}\e^{-2}+\e^{-1}+1\right)-2 \tan^{-1}(\e^{-1}+1)+\frac{\pi}{4}+ \frac{\epsilon}{1+2\e+\e^2} .\]
As $\e\to 0$, the first term goes to infinity and the other terms converge to finite limits. Therefore, the integral $I$ goes to infinity as $\e\to 0$.  On the other hand,
$$
\int \frac{x^2 dx}{((x+\e)^2+\e^2)^2}= \frac{1}{\e}\tan^{-1} \left(\frac{x+\e}{\e}\right)+
\frac{\e^2}{x^2+2x\e+2\e^2},
$$
so that 
\[ II= 2\left(\tan^{-1}(\epsilon^{-1}+1)-\frac{\pi}{4}\right) + 2\e\left(\frac{\epsilon^2}{1+2\e+2\e^2}-\frac{1}{2}\right).\] 
As $\e\to 0+$,  the integral $II$ converges to the limit $\dfrac{\pi}{2}$.
This shows that the limit in~\eqref{nolimit} does not exist,  since its real part goes to $+\infty$ as $\e\to 0$.
Therefore ${\bc}\,f$ cannot be defined. 

 \end{proof}

Consider now the domain of the form $\Omega \times \mathbb C \subset \Mm=\mathbb C^2_{(z_1,z_2)}$ which does not have a generic corner at the origin. From above computations,  it follows that the function $\frac{1}{z_1^2}$ does not admit the boundary current. This gives examples of nonexistence at nongeneric corners in higher dimensions.
\section{An open problem: the global extension phenomenon}
One of the  important aspects of the theory of boundary values is the reconstruction property, i.e., restoring
the function from its values on the boundary.  Such a problem can be posed in both a local and global version.
For a CR function on the smooth connected boundary of a domain in 
$\mathbb C^n$ the global extension to the domain as a holomorphic function  may be obtained by means of the Bochner-Martinelli integral (see, e.g., \cite{kytbook}). This is 
known in the literature as the Bochner-Hartogs phenomenon, and can be viewed as a generalization of classical
Hartogs' Kugelsatz. For boundary currents defined as in Theorem~\ref{thm-bcexistence} the problem is
two-fold: first one needs to identify the class of currents in $\Dd'_{0,1}(\Mm)$ that are boundary values of holomorphic functions of polynomial growth (i.e., to determine the range of the operator $\bc$), and secondly  to reconstruct the holomorphic function given any current in that class. 
While this problem is open for general piecewise smooth domains,
in \cite{chak-sh} we are able to solve it for product domains. Here we give a short account of our result, the details may be found in \cite{chak-sh}.

Let $\Mm_1,\dots,\Mm_N$ be complex manifolds, and $\Mm=\Mm_1\times\dots\times\Mm_N$.
Let $D_j\Subset\Mm_j$ be a domain with $\mathcal{C}^\infty$-smooth boundary, $j=1,\dots, N$. 
Then $\Omega=D_1\times\dots \times D_N$ is a product domain in our sense. We also set
\begin{equation}\label{eq-omegaj}
\Omega_j = \Mm_1\times\dots \times D_j \times\dots\times\Mm_N,
\end{equation}
 and observe that $\Omega = \bigcap_{j=1}^N \Omega_j$.
It is easy to see that each corner is a CR manifold, and so $\Omega$ has generic corners.
We define the subspace $\mathcal{Y}^{0,1}_\Omega(\Mm)$ of
$\Dd'_{0,1}(\Mm)$ as follows. A current $\gamma\in \Dd'_{0,1}(\Mm)$ belongs to $\mathcal{Y}^{0,1}_\Omega(\Mm)$ if the following conditions are satisfied:

\begin{enumerate}\item 
$\gamma$ satisfies the {\it Weinstock condition} with respect to $\Omega$, i.e., 
for $\omega\in \Dd^{n,n-1}(\Mm)$, we have
\begin{equation}\label{eq-weinstock}
\dbar\omega=0 {\rm\ on\ } \overline{\Omega} \ \ \Longrightarrow\ \ \pair{\gamma}{\omega}=0.
\end{equation}
This is a generalization of the usual tangential Cauchy-Riemann equations 
for the boundary values of holomorphic functions, in fact, for domains in $\mathbb C^n$ with connected complement, the Weinstock condition is equivalent to $\gamma$ being $\overline\partial$-closed.

\item Suppose that the piecewise smooth domain $\Omega$ is represented as an intersection of smoothly bounded domains. Let
\begin{equation}\label{eq-iotaj}
\iota^j:\partial\Omega_j\to \Mm, \ \ j=1,\dots, N,
\end{equation}
be the inclusion maps. Then there exist distributions $\alpha_j\in \Dd'_0(\partial\Omega_j)$ with support in $\partial\Omega_j\cap \overline{\Omega}$ such that
we can write
\begin{equation}\label{eq-facewise}
\gamma=\sum_{j=1}^N\left(\iota^j_*(\alpha_j)\right)^{0,1}.
\end{equation}
We will call the distributions $\alpha_1,\dots, \alpha_N$ the {\em face distributions} associated with the 
current~$\gamma$.

\item The third condition, which we call {\it canonicality of face distributions} is rather technical, and cannot
be stated precisely without introducing some relevant technical notions. A full explanation  may be found in \cite{chak-sh}. 
Informally, it can be understood as follows. Given a function $f\in \mathcal{A}^{-\infty}(\Omega)$ on a smooth domain, there exists the extension of $f$ as a distribution in $\mathcal{D}'_0(\Mm)$ with the property that it vanishes outside $\overline \Omega$ and its values on $\partial \Omega$ are determined in a limit process from the values in $\Omega$, similar to that in Theorem~\ref{thm-bcexistence}. This is called the {\it canonical extension} of $f$. A similar canonical extension exists for the distributions 
$\alpha_j\in \Dd'_0(\partial\Omega_j)$ defined by~\eqref{eq-facewise}. The condition now
is that the canonical extensions of $\alpha_j$ agree with $\alpha_j$. In particular, this condition ensures that one can talk about boundary values of the face distributions themselves along higher codimensional strata. 

\end{enumerate}

We note that all three conditions above are satisfied by boundary currents of holomorphic functions. In fact,
we have the following characterization of the distributional boundary values of holomorphic functions on product domains:

\begin{thm}\label{thm-product} Let $\Omega$ be a product domain as above. Then for each $f\in \mathcal{A}^{-\infty}(\Omega)$, we have $\bc f\in \mathcal{Y}^{0,1}_\Omega(\Mm)$, and 
the map
\[ \bc: \mathcal{A}^{-\infty}(\Omega)\to \mathcal{Y}^{0,1}_\Omega(\Mm)\]
is an isomorphism of topological vector spaces.
\end{thm}

We remark that for a smoothly bounded domain $\Omega$ the third condition is void, and the second
condition simply means that there exists a distribution $\alpha\in \Dd'_0(\partial\Omega)$ such that 
$\gamma=\iota_*(\alpha)^{0,1}$. This has a simple geometric interpretation: if a $(n,n-1)$-form $\phi$
vanishes on $\partial\Omega$, then $\gamma(\phi)=0$. In particular this means that for smoothly bounded domains in $\mathbb C^n$, the boundary values of holomorphic functions defined as currents are completely equivalent to boundary values viewed as CR distributions.

Note that conditions (1) and (2) above make sense in any piecewise smooth domain. Therefore, 
we can formulate a more precise version of the problem of global extension in the following form:

\medskip

{\bf Open problem:} {\em Let $\Mm$ be a complex manifold, and let $\Omega\Subset \Mm$ be a domain with generic corners.  Let $\gamma\in \Dd'_{0,1}(\Mm)$ be a current which satisfies the conditions (1) and (2) above, i.e., the Weinstock condition, and the fact that $\gamma$ can be represented in terms of face distributions $\alpha_j$ on the faces of the domain. What further condition do we need to impose on $\gamma$, so that there is a holomorphic function 
 $f$ on $\Omega$ with $\gamma=\bc f$?}



\begin{thebibliography}{BER99}

\bibitem[BER99]{baouendibook}
M.~Salah Baouendi, Peter Ebenfelt, and Linda~Preiss Rothschild.
\newblock {\em Real submanifolds in complex space and their mappings},
  volume~47 of {\em Princeton Mathematical Series}.
\newblock Princeton University Press, Princeton, NJ, 1999.

\bibitem[CS]{chak-sh}
Debraj Chakrabarti and Rasul Shafikov.
\newblock Distributional boundary values of holomorphic functions on product
  domains.(Preprint).
  \newblock Available online at \url{http://arxiv.org/abs/1505.01230}.

\bibitem[HL75]{hala1}
F.~Reese Harvey and H.~Blaine Lawson, Jr.
\newblock On boundaries of complex analytic varieties. {I}.
\newblock {\em Ann. of Math. (2)}, 102(2):223--290, 1975.

\bibitem[Kyt95]{kytbook}
Alexander~M. Kytmanov.
\newblock {\em The {B}ochner-{M}artinelli integral and its applications}.
\newblock Birkh{\"a}user Verlag, Basel, 1995.
\newblock Translated from the Russian by Harold P. Boas and revised by the
  author.

\bibitem[PW78]{polkingwells}
John~C. Polking and R.~O. Wells, Jr.
\newblock Boundary values of {D}olbeault cohomology classes and a generalized
  {B}ochner-{H}artogs theorem.
\newblock {\em Abh. Math. Sem. Univ. Hamburg}, 47:3--24, 1978.
\newblock Special issue dedicated to the seventieth birthday of Erich
  K{{\"a}}hler.

\bibitem[Str84]{straube84}
Emil~J. Straube.
\newblock Harmonic and analytic functions admitting a distribution boundary
  value.
\newblock {\em Ann. Scuola Norm. Sup. Pisa Cl. Sci. (4)}, 11(4):559--591, 1984.

\end{thebibliography}
\end{document}